\documentclass[11pt,draft]{amsart}
\usepackage{amssymb,amsxtra}
\usepackage{verbatim,fullpage,enumerate}
\usepackage[english]{babel}
\usepackage[T1]{fontenc}
\usepackage{xspace}

\numberwithin{equation}{section}

\newtheorem{Theoreme}{Theorem}[section]
\newtheorem{Prop}[Theoreme]{Proposition}
\newtheorem{prop}[Theoreme]{Proposition}

\newtheorem{Lemme}[Theoreme]{Lemma}

\newtheorem{cor}[Theoreme]{Corollary}

\theoremstyle{definition}

\theoremstyle{remark}
\newtheorem{Remarque}[Theoreme]{Remark}

\newenvironment{dem}{\proof[Proof]}{\endproof}

\def\preisomto{\vbox{\hbox to
                 15pt{\hfill$\sim$\hfill}\nointerlineskip\vskip
                 -0.3pt
                 \hbox to 15pt{\rightarrowfill}}}

\def\prelongisomto{\vbox{\hbox to
                18pt{\hfill$\sim$\hfill}\nointerlineskip\vskip -0.3pt
                \hbox to 18pt{\rightarrowfill}}}

\def\longisomto{\mathop{\prelongisomto}}

\def\Ker{\operatorname {Ker}}

\def\Hom{\operatorname {Hom}}
\def\End{\operatorname {End}}
\def\Aut{\operatorname {Aut}}

\def\Int{\operatorname {Int}}
\def\id{\operatorname {id}}
\def\lg{\operatorname {lg}}

\def\GL{\operatorname {GL}}

\newcommand{\ifff}{\Longleftrightarrow}

\newcommand{\N}{{\mathbb{N}}}
\newcommand{\C}{{\mathbb{C}}}

\newcommand{\Z}{{\mathbb{Z}}}
\newcommand{\D}{{\mathbb{D}}}
\renewcommand{\H}{{\mathbb{H}}}

\newcommand{\U}{\mathbb{U}}

\newcommand{\hbi}{\overline{\mathbb{H}_k(\mu_\varpi)}}

\newcommand{\CT}{{\mathbb{C} T}}

\newcommand{\che}{\mathcal{H}}
\newcommand{\calP}{\mathcal{P}}
\newcommand{\calR}{\mathcal{R}}
\newcommand{\calS}{\mathcal{S}}

\newcommand{\calC}{\mathcal{C}}
\newcommand{\calCD}{\mathcal{C}\spcheck}

\newcommand{\calA}{\mathcal{A}}

\newcommand{\calL}{\mathcal{L}}

\newcommand{\tk}{\tilde{k}} 
\newcommand{\tgamma}{{\tilde{\gamma}}}

\newcommand{\tnu}{{\tilde{\nu}}}

\newcommand{\vkappa}{\varkappa}

\newcommand{\twisted}[2]{{}^{{#1}}\!{{#2}}}
\newcommand{\RMod}[1]{{{#1}\text{-}\mathrm{Mod}}}
\newcommand{\Rmod}[1]{{{#1}\text{-}\mathrm{mod}}}
\newcommand{\ModR}[1]{{\mathrm{Mod}\text{-}{#1}}}

\begin{document}

\title{Generalized graded Hecke algebra for complex reflection group of type G$(r,1,n)$}

\author{C. Dezélée}

 \address{LMV, Universit{\'e} de Versailles, France}
 \email{dezelee@math.uvsq.fr}

\date{4 mai 2006}

\subjclass[2000]{16Sxx, 20Cxx, 17Bxx}

\keywords{Complex reflexion groups, Hecke algebras, principal series, symplectic reflection algebra}

\maketitle

\begin{abstract}

We study an algebra similar to a graded Hecke algebra, associated
to a complex reflection group of type G(r,1,n) and isomorphic to a
subalgebra of the symplectic reflection algebra. We define
principal series modules over it and prove an irreducibility
criterion for them.

\end{abstract}


\addtocounter{section}{-1}



\section{Introduction}

\label{Intro}

 Let $G \subset U(V)$ be a finite complex reflection group acting over a $n$-dimensional hermitian space $V$.
 We denote by $\calR$ the set of its reflections. Recall that $G$ acts naturally over $\calP=S(V)$,
 $\calS=S(V^*)$, and over $\calR$ by conjugation. Let $k$~$:$~$\calR \rightarrow \C$ be a $G$-invariant
parameter-function. We can consider two $\C$-algebras associated to these data:\\

$\bullet$ The symplectic reflection algebra $\che(G,k)$, also called rational Cherednik algebra
in this case (see \cite{EG}, \cite{G}, \cite{GGOR}) which is a two step degeneration of the double affine Hecke algebra,
isomorphic to $\calP\otimes\C G \otimes \calS$ as a vector space.\\

$\bullet$ The graded Hecke algebra $\H_{Dr}(G,k)$ defined by Drinfeld (cf.~\cite{Dr}) and isomorphic to $\calP \otimes \C G$ as a vector space.\\

It is well known that $\H_{Dr}(G,k)$ is isomorphic to a subalgebra
of $\che(G,k)$ when $G$ is a Weyl group of type $A$. In
\cite{De2}, we introduced a generalized graded Hecke algebra which
can be realized as a subalgebra of $\che(G,k)$ when $G$ is a Weyl group of type $B$ or $D$,
and studied its representation theory.\\




Moreover, when $G$ is a certain kind of wreath-product, the
representation theory of the Cherednik algebra $\che(G,k)$ is
linked to the ring of polynomial differential operators over the
space of representations of a quiver associated to $G$ via the
McKay correspondance (cf.\cite{EGGO}, \cite{Go}).

This motivates the research of a generalization of graded Hecke algebras associated to a
wreath-product which could be realized in the associated symplectic reflection algebra.\\

In \cite{RS}, A. Ram and A.V. Shepler classified the graded Hecke
algebras associated to complex reflection groups obtained by
Drinfeld's construction. For type $G(r,1,n)$, $r>2$, $n>3$, there
is no non-trivial construction. In these cases, the authors define
a new algebra $\H(r,n)$, isomorphic to $S(V)\otimes \C G$ and
containing the algebra $\C G$.\\

In this paper, we show that the algebra defined by Ram and Shepler
(up to a deformation  by the parameter $k$) can be realized in the
associated symplectic reflection algebra $\che(G,k)$. This allows us to study
representation theory for this new algebra $\H_k(r,n)$, as in \cite{De2}.\\





In section \ref{Prelim}, we recall definition and PBW theorem for
the symplectic reflection algebra and give the definition of the
generalized graded Hecke algebra $\H_k(r,n)$ associated to a
complex reflection groups $G$ of type $G(r,1,n)$. In paragraph
\ref{Real} we realize $\H_k(r,n)$ as a subalgebra of $\che(G,k)$
using a convenient family of commuting elements in $\che(G,k)$.
This allows us to prove a PBW for $\H_k(r,n)$, to calculate its
center $Z_k(r,n)$ and to identify its largest commutative
subalgebra $\calC \simeq \calP^G \otimes Z(\C G)$
(cf.~section~\ref{BP}). Then we can define the principal series
module $M_k(\gamma\otimes \tau)$ by induction from the irreducible
$\calC$-module of character $\gamma\otimes \tau$ ;
each irreducible $\H_k(r,n)$-module is a quotient of a principal series module
(cf. section~\ref{PSMdef}).\\
An irreducibility criterion for $M(\gamma \otimes \tau)$,
involving $k$, $\tau$ and $\gamma$, is given in
Theorem~\ref{thm13}, by reduction to the case of a principal
series module over a graded Hecke algebra associated to a Weyl
group.

\section{Preliminaries}

\label{Prelim}

Let $V$ be a $n$-dimensional hermitian space and $G$ a complex reflection group of type $G(r,1,n)$.
We denote by $\calA$ the set of hyperplanes of complex reflections in $G$. Let
$\xi$ be a primitive root of unity of order $r$.\\
There exists an orthonormal basis $v_1,\dots, v_n$ of $V$ such
that $\calA=\{H_s,\, H_{i,j}^t / \, 1 \le s, i \neq j \le n, \, 0
\le t \le r-1\},$ where  $H_s=Vect\{v_i/\, i \neq s\}$ and
$H_{i,j}^t=Vect\{\xi^t v_i +v_j, \, v_s / \, s \neq i, \, s \neq j \}$.\\
We denote by $\xi_i$ the complex reflection  with $H_i$ as
eigenspace for the eigenvalue $\xi$, and by $s_{u,v}^{[j]}$ the
only complex
reflection stabilizing $H_{u,v}^j$, for $1 \le i, \quad u \neq v \le n$, $0\le j \le r-1$.\\
Notice that the set $\calR$ of complex reflections in $G$ is
$\{\xi_i^t, \, s_{u,v}^{[m]} / \, 1 \le i,u\neq v \le n, \, 1 \le t \le r-1, \, 0 \le m \le r-1\}$.\\
Let $k$ $:$ $\calR : \rightarrow \C$ be a $G$-invariant
parameter-function. For each $t \in \{1,\dots,r-1\}$ we denote by
$k_t$ its value over all $\xi_i^t$, $1 \le i \le n$, and we denote
by $\bar{k}_0$ its value over all $s_{u,v}^{[j]}$, $1\le u \neq v \le n$, $j \in \{0, \dots,r-1\}$.\\
Let $\alpha_1,\dots,\alpha_n \in V^*$ be the dual basis of
$v_1,\dots,v_n$ according to the hermitian structure of $V$. We
note $v_{i,j}^{[t]}=\displaystyle{\frac{\xi^{-t}
v_i-v_j}{\xi^{-2t}+1}}$ and
$\alpha_{i,j}^{[t]}=\xi^{-t}\alpha_i-\alpha_j$, for $1 \le i \neq
j \le n$ and $1 \le t \le r-1$.

\subsection{Symplectic reflection algebra}
\label{SRA}

The symplectic reflection algebra $\che(G,k)$ associated to these
data is generated by $x \in V^*$, $y \in V$ and $g \in G$ with the following
defining relations (cf. \cite{EG}, \cite{G},\cite{GGOR}):\\
$g.x.g^{-1}=g(x),\quad  g.y.g^{-1}=g(y), \quad [x_1,x_2]=0, \quad  [y_1,y_2]=0, \quad  \forall y_1,y_2 \in V, \quad \forall x_1,x_2 \in V^*,$\\
$[y,x]=  <y,x>- \displaystyle{\sum_{t=1}^{r-1}
k_t\sum_{i=1}^{n}<y,\alpha_i><v_i,x>\xi_i^t}
 -\bar{k}_0 \displaystyle{\sum_{m=0}^{r-1} \sum_{1 \le i\neq j \le n}
<y,\alpha_{i,j}^{[m]}> <v_{i,j}^{[m]},x> s_{i,j}^{[m]}}.$\\

This algebra contains $\calP$, $\calS$ and $\C G$ as subalgebras
and verify a PBW property (see e.g. \cite[Theorem~1.3]{EG} ),
which can be expessed by the following isomorphism of vector
spaces:

\begin{equation}
\label{PBW0}
\che(G,k) \simeq \calP\otimes \C G \otimes  \calS.
\end{equation}

For any complex reflection group $G \subset U(V)$ and any
function-parameter $k$, we denote by $\D(G,k)$ the subalgebra of
$\che(G,k)$ generated by $G$ and the elements $v_i \alpha_i$, $1
\le i \le n$.

This algebra is a generalization for complex reflection group of
the algebra $\D$ defined in \cite{De2} for $G$ a Weyl group. It
has been noticed (e.g. in \cite{BEG}) that $\D$ is a graded Hecke
algebra, as defined by Lusztig, when $G$ has type $A$ ($r=1$).

\subsection{Generalized graded Hecke algebra}
\label{GGHA}

Let $V_0$  be a n-dimensional hermitian space and $G \subset
U(V_0)$ a complex reflection group of type $G(r,1,n)$. Then $V_0$
is endowed with $e_1, \dots, e_n$ an orthonormal basis such that
for $u \neq v \in \{1,\dots,n\}$, the reflection $(u,v)$
exchanging $e_u$ and $e_v$ belongs to $G$, and for each $i \in
\{1,\dots n\}$, the reflection $\theta_i$, of matrix
$\rm{diag}(1,\dots,1,\xi,1,\dots,1)$
in the basis $e_1,\dots,e_n$ with $\xi$ in $i^{\rm{th}}$ position, belongs to $G$.\\
We denote by $S_n$ the symmetric group generated by the $(u, v) \in G$,
and by $T \simeq {(\Z/r\Z)}^n$ the subgroup of $G$ generated by the $\theta_i$, $1 \le i \le n$.\\

Notice that the reflections $\theta_1$ and $(i,i+1)$, $1 \le i \le n-1$, generate $G$,
and that $G \simeq T \rtimes S_n.$\\
Fix $c_0 \in \C$.\\
We define a function-parameter $\tilde{c}$ from the set of
reflections of $S_n$ to $\C T$. For each $1 \le u \neq v \le n$,
we note :
\begin{equation}
\label{mulHgr}
\tilde{c}_{(u,v)}=\displaystyle{\frac{c_0}{2}\sum_{i=i}^{r-1}
\theta_u^i \theta_v^{-i}+\theta_u^{-i}\theta_v^i.}
\end{equation}
Notice that $\tilde{c}$ is $S_n$-equivariant and that in the basis
$e_1, \dots, e_n$, $\tilde{c}_{(i,j)}$ has a diagonal matrix with
$0$ in the $i^{\rm{th}}$ and $j^{\rm{th}}$ positions, and with $1$ anywhere else.\\

 The algebra $\H(G,c_0)$ is generated by the $e_1,\dots,e_n$ and $t_g$, $g \in G$,  with the following
relations :

\begin{equation}
\label{Hgr1}
[e_i,e_j]=0, \quad t_{u}t_v=t_{uv}\quad   \forall i,j , \quad \forall u,v \in G
\end{equation}

\begin{equation}
\label{Hgr2}
t_{\theta_i} e_j=e_j t_{\theta_i} \quad   \forall i,j \in \{1,\dots,n\},
\end{equation}

\begin{equation}
\label{Hgr3}
t_{(j,j_+1)} e_i= e_i t_{(j,j+1)}  \quad \forall i,j, / \quad i\neq j, \quad  i\neq j+1,
\end{equation}

\begin{equation}
\label{Hgr4}
t_{(i,i+1)}e_i=e_{i+1}t_{(i,i+1)}-\tilde{c}_{i,i+1}t_{(i,i+1)}, \quad \forall i \in \{1,\dots,n-1\}.
\end{equation}

Notice that when $r=1$, $\H(G,k)$ is the graded Hecke algebra of
type $A$ defined by Lusztig (\cite{Lu}). When $r=2$, $\H(G,k)$ is
the generalization of graded Hecke algebra studied in \cite{De2}.

\subsection{Realization in the symplectic reflection algebra}

\label{Real}

From now on, $G$ is a complex reflection group of type $G(r,1,n)$, we use the previous notations.
For $j \in \{1,\dots,n\}$, we define the following elements in
$\che(G,k)$ :

\begin{equation}
\label{defopka} D_j(G,k)=\displaystyle{v_j \alpha_j - \sum_{1 \le
i <j} \tk_{(i,j)}s_{i,j}^{[0]}},
\end{equation}
where $\tk_{(i,j)}=\displaystyle{ \frac{\bar{k}_0}{2}
\sum_{m=0}^{r-1} \xi_i^m\xi_j^{-m} +\xi_j^m\xi_i^{-m}}$.

We will simply denote $D_i(G,k)$ by $D_i$.
As in the case $r=2$ (cf. \cite{Ka} and \cite{De2}), these
operators verifies the following relations in $\che(k)$ :

\begin{equation}
\label{Dgr1}
D_i D_j=D_j D_i, \quad   \forall i,j \in \{1,\dots,n\}.
\end{equation}

\begin{equation}
\label{Dgr2}
D_i \xi_j=\xi_j D_i, \quad   \forall i,j \in \{1,\dots,n\}.
\end{equation}

\begin{equation}
\label{Dgr3}
s_{j,j_+1}^{[0]} D_i= D_i s_{j,j_+1}^{[0]}  \quad \forall i,j, \quad i\neq j, \quad  i\neq j+1,
\end{equation}

\begin{equation}
\label{Dgr4}
s_{i,i_+1}^{[0]}D_i=D_{i+1}s_{i,i_+1}^{[0]}-\tilde{k}_{i,i+1}s_{i,i_+1}^{[0]}
,\quad\ \forall i \in \{1,\dots,n-1\}.
\end{equation}

Notice that $\D(G,k)$ is the subalgebra of $\che(G,k)$ generated
by $G$ and the elements $D_i$, $1 \le i \le n$. These elements
allow us to realize $\H(G,\bar{k_0})$ in $\che(G,k)$.

\begin{prop}
\label{propreal}
There is an isomorphism of algebras $\psi$ $:$ $\H(G,k) \rightarrow \D(G,\bar{k_0})$ given by $\psi(e_i)=D_i(G,k)$, $\psi(t_{\theta_i})=\xi_i$, $\psi(t_{(u,v)})=s_{u,v}^{[0]}$ for
$1\le i, \quad u \neq v \le n$.
\end{prop}

\begin{dem}
Since the $D_i$'s are clearly linearly independant, we only have
to verify that the $\D_i$'s and the elements of $G$ obey the
relations \ref{Hgr1}, \ref{Hgr2}, \ref{Hgr3}, \ref{Hgr4} in
$\che(k)$, which is a obvious with relations \ref{Dgr1},
\ref{Dgr2},\ref{Dgr3}, \ref{Dgr4}.
\end{dem}

Notice that $\H(G,c_0)$ can be realized
in any algebra $\che(G,k)$ such that $\bar{k_0}=c_0$. From now on, we
will denote $\H(G,\bar{k_0})$ by $\H_k(r,n)$ and identify it to the
subalgebra $\D(G,k)$ of $\che(G,k)$.

\section{Basic properties}

\label{BP}

Thanks to the previous realization, we can prove a PBW property
for $\H_k(r,n)$.

\begin{prop}
\label{PBW1}
As vector-spaces : $\H_k(r,n) \simeq S(V_0) \otimes \C G.$
\end{prop}
\begin{dem}
From now on, we can identify $V_0$ to the subspace
$\bigoplus_{j=1}^n \C D_j$ in $\che(G,k)$, via $\psi$. Since the
$D_i$'s commute and since $\che(G,k)$ verify a PBW property, they
generate a subalgebra of $\che(G,k)$ isomorphic to $S(V_0)$. Then
the proposition is shown in a similar way than in \cite{De2},
Lemma~1.2 and Theorem~1.3.
\end{dem}

Let us precise now the commutation relations in $\H_k(r,n)$  as we
did in \cite[Lemma~1.4]{De2}. For $w \in S_n$and $i,j \in
\{1,\dots,n\}$ , we note $w*i=j$ if $w(e_i)=e_j$.
Recall that the inversion set of $w \in S_n$ is $R(w) = \{\{i,j\},
\, 1 \le i<j \le n / \, w*j < w*i\}$.

\begin{Lemme}
\label{lemcom} For all $\zeta \in V_0$ and $w \in S_n$, we have
the following relation in $\H_k(r,n)$:\\
$ w \zeta = \displaystyle{w(\zeta)w - \sum_{\{i,j\} \in R(w^{-1})}
<v_{i,j}^{[0]},w(\zeta)>(i,j)w\tk_{w^{-1}*i,w^{-1}*j}}$.
\end{Lemme}
\begin{dem}
  Since the multiplicity $\tk$ is $S_n$-equivariant and
commutes with the elements of $V$, one can prove the formula in
the same way as in~\cite[Proposition~1.1(1)]{Op2} in the graded
Hecke algebra case.
\end{dem}

Notice that the center $Z(G)$ of $\C G$ is $\C T^{S_n}$, since the
elements of $T$ are $S_n$-equivariant and since $T$ is the larger
commutative subgroup in $G$.
 Thanks to Proposition~\ref{PBW1}, Lemma~\ref{lemcom} and with the method used in \cite[Theorem~1.8]{De2}, we can now calculate the center of $\H_k(r,n)$.

\begin{Prop}
\label{centerHgr]}
The center $Z_k(r,n)$ of the algebra $\H_k(r,n)$ is isomorphic via
$\psi$ to $S(V_0)^{S_n} \otimes \C T^{S_n}$.
\end{Prop}

From now on we simply denote by $(i,j) \in \che(G,k)$ the
reflection $s_{i,j}^{[0]}$ exchanging $v_i$ and $v_j$ in $V$.\\

\section{Representation theory of $\H_k(r,n)$}
\label{RTH} Generalizing the case of type $A$ ($r=1$), we define
some principal series module over $\H_k(r,n)$.

\subsection{Weights of finite dimensional $\H_k(r,n)$-modules}
\label{weight}

Notice that the larger commutative subalgebra $\calC$ in
$\H_k(r,n)$ is $S(V_0) \otimes \C T$, so any algebra morphism $\tgamma : S(V_0)\otimes \CT\rightarrow \C$ is determined by a linear form $\gamma \in V_0^*$ and a  character $\mu \in T\spcheck$; we will denote it by
$\gamma\otimes\mu$.  Thus we have $(\gamma\otimes\mu)(v \oplus t)=\gamma(v)\mu(t)$, for all $v \in V_0$, $t \in \CT$.
We denote by
$\calCD =\{\gamma\otimes \mu : \gamma \in V_0^*, \quad \mu \in T\spcheck\}$
the set of characters of $\calC$.\\

Let $M$ be a finite dimensional $\H_k(r,n)$-module and fix $\tnu \in \calCD$.
The weight subspaces and generalized weight subspaces of $M$ associated to $\tnu$ are respectively:
$$
M_{\tnu}=\{ m \in M : \forall a \in \calC, \quad a.m=\tnu(a)m \}
$$
and
$$
M_{\tnu}^{\rm{gen}}=\{ m \in M : \forall a \in \calC, \quad \exists k \in \N, (a-\tnu(a))^k.m=0 \}.
$$

We will say that $\tnu$ is a weight of $M$ if $M_{\tnu}^{\rm{gen}}
\ne 0$ (or equivalently $M_\tnu \ne 0$).

\subsection{Automorphisms and anti-automorphisms in $\H_k(r,n)$}
\label{dual} Let us recall some notations and facts used in
\cite{De2}. For any $\C$-algebra $R$, we denote by $\RMod{R}$,
resp.~$\ModR{R}$, the category of the left, resp.~right,
$R$-modules. The lattice (for the inclusion order) of the
submodules of $M \in \RMod{R}$ is denoted by $\calL_R(M)$, and the
length of $M$ is denoted by $\lg_R(M)$.

Let $S$ be another $\C$-algebra and $N \in \RMod{S}$; if there
exists an isomorphism of $\C$-vector spaces $f: M \to N$
such that $X \mapsto f(X)$ is a lattice isomorphism between
$\calL_R(M)$ and $\calL_S(N)$, then we will denote this fact by
$\calL_R(M) \cong \calL_S(N)$.\\

If $M,N \in \RMod{R}$ have finite length and if $\calL_R(M) \cong
\calL_S(N)$ by $f$, then each composition series of $M$ is clearly
sent by $f$ onto a composition series of $N$, thus the $\lg_R(M)=\lg_R(N)$.\\

Recall that for all $M \in \RMod{R}$ and $\vkappa \in \Aut_\C(R)$,
we can define a twisted module $\twisted{\vkappa}{M} \in \RMod{R}$
in the following way: $\twisted{\vkappa}{M} = M$ as a
$\C$-vector space, endowed with the action $a\centerdot u
=\vkappa(a).u$ for $a \in R, u \in M$. Observe that any
anti-automorphism $\iota$ on the $\C$-algebra $R$
provides an isomorphism between $\RMod{R}$ and $\ModR{R}$: we turn
$M \in \RMod{R}$ into a right $R$-module $M^{\iota}$
by puting $M^{\iota} = M$ as $\C$-vector space and $u.a =
\iota(a).u$ for all $u \in M$, $a \in R$.\\
Then the dual $M^* = \Hom_\C(M,\C)$ can be endowed with a left
$R$-module structure by the formula $ < a.f, u > \quad = \quad < f, \iota(a).u > $
for all $a \in R, f \in M^*, u \in M$.\\

Notice that if $\iota$ is involutive then, for all $M \in\Rmod{R}$, the canonical  bijection
$c_M : M \to M^{**}$ is an $R$-module isomorphism; in particular $\lg_R(M) = \lg_R(M^*)$.\\

These remarks apply to $R=\H_k(r,n)$ in the following cases:

\noindent{(1)} Each $w \in G$ defines an inner automorphism on
$\Int(w) : a \mapsto w aw^{-1}$; we will denote by $\twisted{w}{M}$
the module twisted by $\Int(w)$.

\noindent{(2)} The determinant in $\GL(V_O)$ defines a character
$\det : G \to \U_r \times \{\pm 1\}$, where $\U_r$ denotes the
subgroup of $r^{\rm{th}}$-roots of unity in $\C^*$ . In
particular, $\det$ provides an element of $T\spcheck$.

The map $\delta : x \mapsto \det(x) x$, $x \in T$, induces an
automorphism on $\CT$ such that $\delta(\tk_{i,j})=\tk_{i,j}$ for
all $i,j$. We obtain an automorphism on $\H_k{r,n}$ by extending
$\delta$ as follows:
$$
\delta(v) = v, \quad \delta(\theta) =\det(\theta)\theta, \quad
\delta((i,j)) = (i,j),
$$
for all $v \in V_0,\quad \theta \in T, \quad 1 \le i \neq j \le n$.

\noindent{(3)} One easily checks that the formulas below define an
involutive  anti-automorphism $\iota$ over $\H_k(r,n)$:

\begin{equation}
\label{antiinv}
\iota(D_i) = -D_i, \quad \iota(g) = \det(g)g^{-1},
\quad\quad \forall i \in \{1,\dots,n\}, \quad \forall g \in G.
\end{equation}

In particular, it follows from (3) that $\RMod{\H_k(r,n)} \simeq
\ModR{\H_k(r,n)}$.


\subsection{Principal series modules}
\label{PSM}

We define for $\H_k(r,n)$, as for a graded Hecke algebra, some principal series modules for which we will prove
an irreducibility criterion in the next section.

For all $\tgamma =\gamma \otimes \mu \in \calCD$ we denote by
$\twisted{w}{\tgamma} = \twisted{w}{\gamma} \otimes
\twisted{w}{\mu}$ the diagonal action of $w \in S_n$, and let $\C
m_\tgamma$ be the one dimensional $\calC$-module defined by:
$a.m_\tgamma = \tgamma(a)m_\tgamma$ for all $a \in \calC$. The
principal series $\H_k(r,n)$-module  $M(\tgamma)$ associated to
$\tgamma$ is the induced module:
\begin{equation}
\label{PSMdef}
M(\tgamma)=\displaystyle{\H_k(r,n)\otimes_{\calC}\C m_\tgamma}.
\end{equation}

Clearly, $\{w \otimes m_{\tgamma} : w \in S_n\}$ is a basis of
$M(\tgamma)$, and this basis can be endowed with an order
compatible with the length $\ell$ on $S_n$ defined  by the set of
simple reflections $(i,i+1)$, $1 \le i \le n-1$. According to
Lemma~\ref{lemcom}, for all $\zeta \in V$, $\theta \in T$ and $w
\in S_n$, one has:
\begin{equation*}
\label{eqbase} \zeta \theta w \otimes m_{\tgamma} =
\displaystyle{\twisted{w}{\tgamma}(\zeta \theta) w \otimes
m_{\tgamma}-\sum_{\{u,v\} \in R(w^{-1})}
\twisted{w}{\mu}(\tk_{u,v}\theta)<\alpha\spcheck,\zeta>(u,v)w
\otimes m_{\tgamma}}.
\end{equation*}

Since $\ell((u,v)w)<\ell(w)$ for all $\{u,v\} \in R(w^{-1})$, the
elements of $\calC$ all have an upper triangular matrix in the
basis $\{w \otimes m_{\tgamma} : w \in S_n\}$ ordered by the
length $\ell$, and the weights of $M(\tgamma)$ are the
$\{\twisted{w}{\tgamma} : w \in S_n\}$.\\
The next proposition implies in particular that simple
$\H_k(r,n)$-modules are $n$!-dimensional at most.

\begin{Prop}
\label{propsimple}
Let $M$ be an irreducible $\H_k(r,n)$-module and let $\tgamma$ be a weight of $M$.
Then $M$ is a quotient of $M(\tgamma)$.
\end{Prop}

\begin{dem}
Fix $v_{\tgamma} \in M_{\tgamma} \smallsetminus \{0\}$;
thus $\C v_{\tgamma}$ is an irreducible $\calC$-module of character $\tgamma$. Since the induction is the adjoint functor of the restriction, there exists a unique morphism of $\H_k(r,n)$-modules from $M(\tgamma)$ to $M$ sending $1 \otimes m_\tgamma$ onto $v_\tgamma$.
Since $M$ is irreducible this morphism is surjective and $M$.
\end{dem}

Recall that we can endow the dual of $M \in \Rmod{\H_k(r,n)}$ with a $\H_k(r,n)$-module
structure defined by the anti-automorphism $\iota$, and that we can twist $M$ by $\delta$, see~\ref{dual}.

We will denote by $w_0$ the longest element in $S_n$. For all $\tgamma=\gamma \otimes \mu
\in \calCD$ we put
$$
\tgamma^* = (-\twisted{w_0}{\gamma}) \otimes (-\twisted{w_0}{\mu}).
$$

\begin{Prop}
\label{prop6}
Fix $\tgamma=\gamma \otimes \mu \in \calCD$. Then:
$$
\twisted{\delta}{M(\tgamma)} \simeq M(\gamma \otimes
(-\mu)), \qquad M(\tgamma)^* \simeq M(\tgamma^*).
$$
Each irreducible $\H_k(r,n)$-module is a submodule of a principal series module.
\end{Prop}

\begin{dem}
The proof is similar than in case $r=2$
(cf.~\cite[Proposition~2.2]{De2}).
\end{dem}

\begin{prop}

\label{prop61} Fix $\tgamma=\gamma \otimes \mu \in \calCD$ and $w
\in S_n$. There exists an isomorphism of $\H_k(r,n)$-modules\\
$M(\tgamma) \longisomto \twisted{w}{M(\tgamma)}, \quad t_g \otimes
m_\tgamma \mapsto t_{wg} \otimes m_\tgamma. $
\end{prop}

\begin{dem}
Thanks to Lemma~\ref{lemcom}, the proof is similar than in
\cite[Proposition~2.3]{De2}.
\end{dem}

\subsection{Criterion in Weyl groups' case}
\label{KR}

 Let $V_1$ be a n-dimensional hermitian space. Let $W
\subset U(V_1)$ be a Weyl group of root system $R$ and fix a set
$S$ of simple root in $R$.
Let $c : R \to \C$ a $W$-invariant function. \\

 We recall here the definition and irreducibility criterion of principal series
modules over the graded Hecke algebra associated to $W \subset
U(V_1)$, $S$ and $c$ (cf.~\cite{KR}). We will denote by
$\H_{\mathit{gr}}=\H_{\mathit{gr}}(c,S,V_1)$ this algebra.

We denote by $\mathrm{M}(\lambda)$ the principal series
$\H_{\mathit{gr}}$-module associated to $\lambda \in V_1^*$. We
set
\begin{center}
$P(\lambda)= \{\alpha \in R \quad : \quad \lambda(\alpha)= \pm
c_{\alpha} \}.$
\end{center}

Thus we have the following irreducibility criterion
(cf.~\cite[Theorem~2.10]{KR}): the principal series
$\H_{\mathit{gr}}(c)$-module $\mathrm{M}(\lambda)$ is irreducible
if and only if $P(\lambda)=\emptyset$.

For all $w \in W$ we can twist $\mathrm{M}(\lambda)$ by $\Int(w)
\in \Aut(\H_{\mathit{gr}}(c))$ to define the
$\H_{\mathit{gr}}(c)$-module $\twisted{w}{M(\lambda})$.

\begin{Remarque}

\label{remHgr}

Let $\lambda,\gamma \in V_1^*$. Assume that $\mathrm{M}(\lambda)
\simeq \mathrm{M}(\gamma)$; then $\gamma$ is a weight of
$\mathrm{M}(\lambda)$ and therefore $\gamma \in
W.\lambda=\{\twisted{w}{\lambda} \quad : \quad w \in W\}$. It
follows from~\cite[Proposition~2.7c]{KR} that if
$\mathrm{M}(\lambda)$ is irreducible, so is
$M(\twisted{w}{\lambda})$ for all $w \in W$. Thus, when
$\mathrm{M}(\lambda)$ is irreducible:
$$\mathrm{M}(\lambda) \simeq \mathrm{M}(\gamma) \ifff  \gamma \in W.\lambda.$$

\end{Remarque}

A proof similar to the Proposition~\ref{prop61} yields :

\begin{prop}

\label{prop62}

Fix $\lambda \in V_1^*$ and $w \in W$. There exists a
$\H_{\mathit{gr}}(c)$-module isomorphism
\begin{equation*}
 \mathrm{M}(\lambda)
\longisomto \twisted{w}{\mathrm{M}(\lambda)}, \quad t_g \otimes
m_\lambda \mapsto t_{wg} \otimes m_\lambda.
\end{equation*}
\end{prop}



\section{Irreducibility criterion}

\label{secIC}


\subsection{Notation}

\label{sec3.1}

For any $\mu \in T\spcheck$, we and denote by $S_n(\mu)$ the
stabilizer of $\mu$ in $S_n$. We denote by $\H_k(\mu)$ the
subalgebra of $\H_k(r,n)$ generated by $V_0$, $S_n(\mu)$ and
$T$.\\

Fix $\tnu=\nu \otimes \varpi \in \calCD$.\\

Let $\Psi : \H_k(r,n) \rightarrow \End(M(\tnu))$ be the
representation of $\H_k(r,n)$ in $M(\tnu)$. We denote by $a.v =
\Psi(a)(v)$ the action of $a \in \H_k(r,n)$ on $v \in M(\tnu)$.\\

For each $s \in \{1,\dots,n\}$, we denote by $\varpi_s$ the unique
integer in $\{0,\dots,r-1\}$ such that
$\varpi(\xi_s)=\xi^{\varpi_s}$, and for each $j \in \{0,\dots,
r-1\}$ we denote by $\varpi^{[j]}$ the cardinal of $\{m/ \,
\varpi_m=j\}$. Then $\varpi$ is determined by the $n$-index
$(\varpi_1,\varpi_2, \dots, \varpi_n)$, and its orbit over the
action of $S_n$ is determined by the $r$-index
$\underline{\varpi}=(\varpi^{[0]},\varpi^{[1]}, \dots, \varpi^{[r-1]})$.\\

The $r$-compositions $\underline{m}=(m_0,\dots,m_{r-1})$ of $n$
(with $0 \le m_i\le r-1$) parameterize the set of $S_n$-orbits in
$T\spcheck$.\\

We denote by $\mu_{\underline{m}} \in T\spcheck$ the element of
the $S_n$-orbit corresponding to $\underline{m}$ such that :
$\mu_{\underline{m}}(\xi_i)=1$ for $1 \le i \le m_0$, and
$\mu_{\underline{m}}(\xi_i)=\xi^j$, for
$\displaystyle{\sum_{s=0}^{j-1} m_s} <i \le
\displaystyle{\sum_{s=0}^{j} m_s}$ and $j \in \{1, \dots, r-1\}$.

Notice that for $\varpi \in T\spcheck$, there exists $\sigma \in
S_n$ such that $\twisted{\sigma}{\varpi}$ is the character
determined by the $n$-index $(0^{\varpi^{[0]}},1^{\varpi^{[1]}},
\dots, (r-1)^{\varpi^{[r-1]}})$, that is
$\mu_{\underline{\varpi^{.}}}$. Moreover we can choose $\sigma$
such that $R(\sigma^{-1})$ only contains some elements $\{i,j\}$
with $\varpi_i \neq \varpi_j$.\\

Notice we have, for any $\{i,j\}\in R(\sigma^{-1})$:
\begin{equation}
\label{simplif1}
\varpi(\tk_{i,j})=\displaystyle{\frac{\bar{k}_0}{2}\sum_{m=0}^{r-1}
\xi^{m(\varpi_i-\varpi_j)}+\xi^{m(\varpi_j-\varpi_i)}=0}.
\end{equation}
Conversely For any $(u,v) \in S_n(\varpi)$, since
$\varpi_u=\varpi_v$, we have:
\begin{equation}
\label{simplif2}
\varpi(\tk_{u,v})=\displaystyle{\frac{\bar{k}_0}{2}\sum_{m=0}^{r-1}
\xi^{0}+\xi^{0}=r \bar{k}_0}.
\end{equation}

From now on, we will simply denote
$\twisted{\sigma}{\varpi}=(0^{\varpi^{[0]}},1^{\varpi^{[1]}},
\dots, (r-1)^{\varpi^{[r-1]}})$ by $\mu_{\varpi}$ (which only
depends on the $S_n$-orbit of $\varpi$).

 We denote by $\H_k(\varpi)$ the subalgebra of $\H_k(r,n)$
generated by $V_0$, $S_n(\varpi)$ and $T$. Let
$\{w_1=\id,w_2,..,w_s\}$ be a set of representatives of
$S_n/S_n(\varpi)$.

\begin{Remarque}

\label{rem7}

The characters $\twisted{w_j}{\varpi}$, $j=1,\dots,s$ are pairwise
distinct. A standard argument implies then that the linear forms
$\twisted{w_j}{\mu} \in (\CT)^*$ are linearly independant. Thus
there exists some elements $y_i \in T$, $i=1,\dots,s$, such that
the matrix $[\twisted{w_j}{\varpi}({y_i})]_{1 \le i,j \le s}$ is
invertible.
\end{Remarque}

For all $j \in \{1,..,s\}$, we define the following subspace of
$M(\tnu)$:
\begin{equation*}
E_j(\tnu)=\bigoplus_{w \in S_n(\varpi)} \C {w_j}w\otimes m_\tnu.
\end{equation*}

Notice that
\begin{equation*}
E_1(\tnu)= \bigoplus_{w \in S_n(\varpi)} \C w\otimes m_\tnu ,
\qquad E_j(\tnu)={w_j}.E_1(\tnu).
\end{equation*}

One has $M(\tnu)=\bigoplus_{j=1}^s E_j(\tnu)$ and this
decomposition is the isotypical decomposition of the $T$-module
$M(\tgamma)$: the group $T$ acts on $E_j(\tgamma)$ by the
character $\twisted{w_j}{\varpi}$.\\

Notice that $S_n(\mu_\varpi)$ is a Weyl group $W_\varpi$ of type
$\mathbf{A}_{\varpi^{[0]}-1}\times
\mathbf{A}_{\varpi^{[1]}-1}\times \cdots \times
\mathbf{A}_{\varpi^{[r-1]}-1}$, with the convention that
$\mathbf{A}_p$ does not appear if $p<0$.\\
For each $i \in \{0,\dots,r-1\}$, we put
$p_i=\displaystyle{\sum_{t=0}^{i}\varpi^{[t]}}$
 We denote by
$S=\{\beta_i=v_i-v_{i+1}, \,1 \le i \le n-1\}$ the set of simple
roots in $S_n$ associated to the lenght $\ell$. Then $S_\varpi =S
\smallsetminus \{ \beta_{p_i}, \, 0\le i \le r-1\}$ is a set of
simple roots for $W_\varpi$. Therefore $S_n(\varpi)=\sigma
S_n(\mu_\varpi)\sigma^{1}$ is a Weyl group with a basis of simple
roots given by $\sigma(\S_\varpi)$.



\subsection{Restriction to $E_1(\tgamma)$}

\label{sec3.2}

Recall that for $M \in \RMod{\H_k(r,n)}$ and $w \in S_n$,
$\twisted{w}{M}$ denotes the $\H_k(r,n)$-module twisted by the
action of $\Int(w)$ (see~\ref{dual}).

\begin{Prop}

\label{prop4} Let $\tgamma=\gamma \otimes \varpi \in \calCD$.
 Fix $w \in S_n$. Suppose that $\mu(\tk_{i,j})=0$
for all $\{i,j\} \in R(w^{-1})$.  Then the map $\phi :
M(\twisted{w^{-1}}{\tgamma}) \rightarrow \twisted{w}{M(\tgamma)}$,
given by $\phi(g\otimes
m_{\twisted{w^{-1}}{\tgamma}})={wgw^{-1}}\otimes m_\tgamma$ for
all $g \in S_n$, is an isomorphism of $\H_k(r,n)$-modules.

\end{Prop}

\begin{dem}

  Since $\{g \otimes m_{\tgamma}={w^{-1}gw}
  \centerdot (1 \otimes m_{\tgamma}, \,  g \in W\}$ is a basis of the space $M(\tgamma)=\twisted{w} M(\tgamma)$,
  the element $1 \otimes m_\tgamma$ generates $\twisted{w}{M(\tgamma)}$.  The group $T$ acts on $1 \otimes
  m_\tgamma \in \twisted{w}M(\tgamma)$ by the character $\twisted{w^{-1}}{\varpi}$.
Then using Lemma~\ref{lemcom} and \eqref{simplif1} in this case,
we notice that $1 \otimes m_\tgamma$ has weight
$\twisted{w^{-1}}{\tgamma}$ in  $\twisted{w}{M(\tgamma)}$. By the
universal property of induction, there exists a unique surjective
$\H_k(r,n)$-morphism $\phi$ from $M(\twisted{w^{-1}}{\tgamma})$
onto $\twisted{w}{M(\tgamma)}$ sending $1\otimes
m_{\twisted{w^{-1}}{\tgamma}}$ to $1 \otimes m_\tgamma$. Since
these two modules are both $n!$-dimensional, $\phi$ is bijective.
Now for $g \in S_n$ one has : \\
$\phi(g\otimes m_{\twisted{w^{-1}}{\tgamma}})=\phi(g.(1 \otimes
m_{\twisted{w^{-1}}{\tgamma}})) =g\centerdot \phi(1 \otimes
m_{\twisted{w^{-1}}{\tgamma}})={wgw^{-1}}.(1 \otimes
m_{\tgamma})={wgw^{-1}} \otimes m_{\tgamma},$ which ends the
proof.

\end{dem}

\begin{cor}
\label{cor5}

{\rm (1)} The map $\phi : M(\twisted{\sigma}{\tgamma}) \rightarrow
\twisted{\sigma^{-1}}{M(\tgamma)}$ given by $\phi(g\otimes
m_{\twisted{\sigma}{\tgamma}})={\sigma^{-1}g\sigma}\otimes
m_\tgamma$, for all $g \in S_n$, is an isomorphism of
$\H_k(r,n)$-modules.

\noindent {\rm (2)} There exists an isomorphism:\\
$\twisted{\sigma}{\phi} : M(\twisted{\sigma}{\tgamma}) \longisomto
M(\tgamma)$ such that $\twisted{\sigma}{\phi}(g \otimes
m_{\twisted{\sigma}{\tgamma}}) = {g \sigma} \otimes m_\tgamma$.

\end{cor}

\begin{dem}

  (1) follows directly from~\eqref{simplif1} and previous proposition applied with $w=\sigma^{-1}$.
\noindent (2) Proposition~\ref{prop61} provides an isomorphism
from $\twisted{\sigma^{-1}}{M(\tgamma)}$ on $M(\tgamma)$ which
sends $u \otimes m_\tgamma$ onto ${\sigma u} \otimes m_\tgamma$;
by composing it with $\phi$ we obtain the required  isomorphism
$\twisted{\sigma}{\phi}$.
\end{dem}

The second assertion of Corollary~\ref{cor5} shows that it is
equivalent to study the $\H_k(r,n)$-module $M(\tnu)$ or the module
$M(\twisted{\sigma}{\tnu})$. We will see that $E_1(\tnu)$ can be
endowed with a structure of $\H_k(\varpi)$-module which determines
the structure of $M(\tnu)$ (cf.~Proposition~\ref{prop12}).

\begin{Lemme}
\label{lem6}

The restriction of $\Psi$ to $\H_k(\varpi)$ equips $E_1(\tnu)$
with a structure of $\H_B(\tnu)$-module.

\end{Lemme}

\begin{dem}

Obviously $E_1(\tnu)$ is stable under the actions of $T$ and
$S_n(\varpi)$. Thus it remains to prove that $E_1(\tnu)$ is stable
under the action of $V$; since the elements of $V$ and $T$ commute
in $\H_k(r,n)$, this result follows from the fact that $E_1(\tnu)$
is the isotypical component of type  $\varpi$ of the $T$-module
$M(\tnu)$.














\end{dem}


\subsection{Isomorphism of lattices}

\label{sec3.2bis}

Fix $\tgamma=\gamma \otimes \varpi \in \calCD$.

We want to show that the lattices $\calL_{\H_k(r,n)}(M(\tgamma))$
and $\calL_{\H_k(\varpi)}(E_1(\tgamma))$ are isomorphic. Observe
first that PBW-property for $\H_k(r,n)$ and $\C S_n(\varpi)
\otimes \CT \otimes S(V_0) \subset \H_k(\varpi)$ imply
%

Therefore
\begin{equation}
\label{eq5} \H_k(r,n) =\displaystyle{\sum_{j=1}^s
{w_j}\H_k(\varpi)}.
\end{equation}

\begin{Prop}
\label{prop12} The module $M(\tgamma)$ is isomorphic to $\H_k(r,n)
\otimes_{\H_k(\varpi)} E_1(\tgamma)$. \\

The maps $Y \to \H_k(r,n)\otimes_{\H_k(\varpi)} Y$ and $X \to X
\cap E_1(\tgamma)$ are inverse bijections from
$\calL_{\H_k(\varpi)}(E_1(\tgamma))$ onto $\calL_{\H_k(r,n)}(M(\tgamma))$.\\
In particular, one has:
\begin{center}
$\lg_{\H_k(r,n)} M({\tgamma}) =\lg_{\H_k(\varpi)} E_1(\tgamma).$
\end{center}
\end{Prop}

\begin{dem}
The proof is  similar to the case $r=2$
(cf.~\cite[Proposition~3.5]{De2}).

\end{dem}

\begin{Remarque}

\label{remeq}

Observe that the previous proposition ensures that the
$\H_k(r,n)$-module $M(\tgamma)$ is simple if and only if the
$\H_k(\varpi)$-module $E_1(\tgamma)$ is simple. Using
Corollary~\ref{cor5}, we conclude that the studies of
$M(\tgamma)$, $E_1(\tgamma)$, $M(\twisted{\sigma}{\tgamma})$ and
$E_1(\twisted{\sigma}{\tgamma})$ are all equivalent to each other.
\end{Remarque}


\subsection{Case $\mu=\mu_\varpi$}

\label{sec3.3}

In this section we assume that $\tgamma = \gamma \otimes
\mu_\varpi \in \calCD$.

Denote by $\H_{\mathit{gr}}(\varpi)$ the graded Hecke algebra
associated to the Weyl group $W=S_n(\mu_\varpi)$, to the set of
simple roots $S_\varpi$ in  $V_1=V_0$ and to the constant
multiplicity $c=r\bar{k}_0$ (cf.~\ref{KR}).

For $\gamma \in V_0^*$, $\mathrm{M}_\varpi(\gamma)$ will denote
the principal series module on $\H_{\mathit{gr}}(\varpi)$
associated to $\gamma$ (i.e. $\mathrm{M}_\varpi(\gamma)$ is the
module $\mathrm{M}(\gamma)$ defined in~\cite{KR}.

Let $I=\Ker(\mu_\varpi)$ be the ideal of $\C T$ generated by the
$x-\mu_\varpi(x)$, $x\in T$. Since the elements of $V$ and $T$
commute, the left ideal $\H_k(\mu_\varpi) I$ is a two-sided ideal.

Thus we can define the algebra
\begin{center}
$\hbi=\displaystyle{\H_k(\mu_\varpi) \otimes_{\C T} (\CT/I) =
\H_k(\mu_\varpi)/\H_k(\mu_\varpi)I.}$
\end{center}
Since the group $S_n(\mu_\varpi)$ is generated by the simple
reflections  in $S(\varpi)$, it follows from~\eqref{Hgr4} and
PBW-property  that $\H_k(\mu_\varpi)=S(V_0)\otimes \C
S_n(\mu_\varpi) \otimes \C T$. As $\C T/I$ is a one dimensional
$T$-module (of character $\mu_\varpi$), there exists an
isomorphism of $S(V_0)$-modules:\\
$\hbi \simeq S(V_0) \otimes \C S_n(\mu_\varpi)$.

Notice that the defining relation \eqref{Hgr4} of $\H_k(r,n)$ and
\eqref{simplif2} yield, in $\hbi$:
\begin{equation*}
(i,i+1)\bar{D_i}=\bar{D_{i+1}}(i, i+1)-r\bar{k}_0,
\end{equation*}
for all $i \in \{1,\dots,n\} \smallsetminus \{p_i, \, 0 \le i \le
r-1\}$. (Here $\bar{D_i}$ stands for the image of $D_i$ in the
quotient $\hbi$.)

These relations coincide with those defining the graded Hecke
algebra $\H_{\mathit{gr}}(\varpi)$. Thus we have a surjective
morphism of algebras \\
 $F: \H_{\mathit{gr}}(\mu_\varpi)\rightarrow \hbi$,
defined by $F(\zeta)=\zeta \pmod{\H_k(\mu_\varpi)I}$ and $F(w)=w
\pmod{\H_B(\mu_i)I}$ for $\zeta \in V$, $w \in S_n(\mu_\varpi)$.
Since $\H_{\mathit{gr}}(\varpi)$ and $\hbi$ are both isomorphic to
$S(V_0) \otimes \C S_n(\mu_\varpi)$ as $S(V_0)$-modules, $F$ is an
isomorphism.

Observe that $T$ acts on the $\H_k(\mu_\varpi)$-module $E_1(\tnu)$
by the character $\mu_\varpi$, thus  $E_1(\tgamma)$ can be viewed
as a $\hbi$-module.

\begin{Lemme}

\label{lem9}

The map $f:\mathrm{M}(\gamma)\rightarrow E_1(\tgamma)$ defined by
$f(t_w \otimes m_\gamma)=t_w \otimes m_\tgamma$ is an isomorphism
which intertwines, via $F$, the actions of
$\H_{\mathit{gr}}(\varpi)$ and of $\hbi$.
\end{Lemme}

\begin{dem}

 Thanks to the isomorphism $F$, $E_1(\tnu)$ can be endowed with a  $\H_{\mathit{gr}}(\varpi)$-module
structure for which it is generated by the vector $1 \otimes
m_\tgamma$. This vector has weight $\gamma$ under the action of
$S(V_0)$. Thus the universal property of $\mathrm{M}(\gamma)$
ensures the existence and the surjectivity of the intertwining
operator $f$. Since both $E_1(\tgamma)$ and $\mathrm{M}(\nu)$ are
$|S_n(\mu_\varpi)|$-dimensional, $f$ is bijective.

\end{dem}

By Lemma~\ref{lem9}, we can identify the
$\H_{\mathit{gr}}(\varpi)$-module $\mathrm{M}(\gamma)$ and the
$\H_B(\mu_\varpi)$-module $E_1(\tgamma)$, in particular we have:
\begin{equation}
\label{eq4} \calL_{\H_{\mathit{gr}}(\varpi)}(\mathrm{M}(\gamma))
\cong \calL_{\H_k(\mu_\varpi)}(E_1(\tgamma)), \quad
\lg_{\H_{\mathit{gr}}(\varpi)}\mathrm{M}(\gamma) =
\lg_{\H_k(\mu_\varpi)} E_1(\tgamma).
\end{equation}

Now we apply this to $\tnu=\nu\otimes \varpi \in \calCD$.\\
Applying~\eqref{eq4} to $\tgamma = \twisted{\sigma}{\tnu} =
\twisted{\sigma}{\nu}\otimes \mu_\varpi$ and using
Corollary~\ref{cor5}, one obtains:
\begin{Prop}
 \label{prop11}
 (Same notations.)
One has\\
$ \calL_{\H_k(\varpi)}(E_1(\tnu)) \cong
 \calL_{\H_{\mathit{gr}}(\varpi)}(\mathrm{M}(\twisted{\sigma}{\nu})),
 \quad \lg_{\H_k(\varpi)}E_1(\tnu)=
 \lg_{\H_{\mathit{gr}}(\varpi)}\mathrm{M}(\twisted{\sigma}{\nu}).$\\

In particular, $E_1(\tnu)$ is an irreducible $\H_k(\varpi)$-module
if, and only if, $\mathrm{M}(\twisted{\sigma}{\nu})$ is an
irreducible $\H_{\mathit{gr}}(\varpi)$-module.
\end{Prop}

\subsection{Irreducibility criterion}

\label{sec3.4}

We can now obtain an irreducibility criterion for the module
$M(\tnu)$.

Recall that $\tnu = \nu \otimes \varpi \in \calCD$ and that
$\sigma \in S_n$ denote the specific permutation introduced
in~\ref{sec3.1}. As in~\cite[Theorem~2.10]{KR}, we set, for all
$\gamma \in V_0^*$:

\begin{equation*}
P_\varpi(\gamma) =\{(i,j) \in S_n(\mu_\varpi) : \nu(D_i-D_j)= \pm
r\bar{k}_0\}.
\end{equation*}

\begin{Theoreme}

\label{thm13}

The following assertions are equivalent:

\begin{enumerate}[{\rm (i)}]

\item $M(\tnu)$ is a simple $\H_k(r,n)$-module;

\item $M(\twisted{w}{\tnu})$ is a simple $\H_k(r,n)$-module for
all $w \in S_n$;

\item $P_\varpi(\twisted{\sigma}{\nu}) = \emptyset$.

\end{enumerate}

\end{Theoreme}

\begin{dem}

We have shown that\\

$M(\tnu) \text{ simple } \ifff M(\twisted{\sigma}{\tnu}) \text{
simple } \ifff E_1(\tnu) \text{ simple } \ifff
\mathrm{M}(\twisted{\sigma}{\tnu}) \text{ simple,} $\\

cf.~Remarque~\ref{remeq},~Proposition~\ref{prop11}.

Therefore, by~\cite[Theorem~2.10]{KR}, (i) is equivalent to (iii).

Since (ii) $\Rightarrow$ (i) is obvious, it remains to prove (i) $\Rightarrow$ (ii).

Assume that $M(\tnu)$ is simple, i.e.
$\mathrm{M}(\twisted{\sigma}{\gamma})$ simple. Let $\tau \in S_n$
be such that $\twisted{\tau w}{\mu}= \mu_\varpi =
\twisted{\sigma}{\mu}$. Then $x= \tau w \sigma^{-1} \in
S_n(\mu_\varpi)$ and we know that:\\
$\mathrm{M}(\twisted{w}{\nu}) \text{ simple } \ifff
\mathrm{M}(\twisted{\tau w}{\nu}) \text{ simple }$\\
(cf.~\cite[Proposition~2.7c]{KR}).

Since $x \in S_n(\mu_\varpi))$, it follows from Remark~\ref{remeq}
, that $\mathrm{M}_i(\twisted{\sigma}{\gamma})$ simple implies
$\mathrm{M}(\twisted{x \sigma}{\gamma})=\mathrm{M}(\twisted{\tau
w}{\gamma})$ simple. Hence the result.

 \end{dem}



\bigskip


\end{document}